\documentclass{birkmult}
\usepackage[T1]{fontenc}
\usepackage{amsmath,amsthm,amscd,amssymb}
\usepackage{graphicx}
\usepackage[ruled,lined,linesnumbered]{algorithm2e}

\newtheorem{theorem}{Theorem}[section]

\newtheorem{lemma}[theorem]{Lemma}
\newcommand{\loc}{\textrm{loc}}
\begin{document}
%opening
\title[Convergence of parametrised non-resonant invariant manifolds]{A note on the convergence of parametrised \\ non-resonant invariant manifolds}
%\date{June 15, 2009}

\author{Tomas Johnson}
\address{Department of Mathematics, Uppsala University, Box 480, 751 06 Uppsala, Sweden}
\email{tomas.johnson@math.uu.se}
\author{Warwick Tucker}
\address{Department of Mathematics, Uppsala University, Box 480, 751 06 Uppsala, Sweden}
\email{warwick.tucker@math.uu.se}

\begin{abstract}
Truncated Taylor series representations of invariant manifolds are abundant in numerical computations. We present an aposteriori method to compute the convergence radii and error estimates of analytic parametrisations of non-resonant local invariant manifolds of a saddle of an analytic vector field, from such a truncated series. This enables us to obtain local enclosures, as well as existence results, for the invariant manifolds.
\end{abstract}


\maketitle

%------------------------------------------------------------------------------------------------------------------------------------------------

\section{Introduction}
The invariant manifolds of a saddle of a vector field are very important objects for the understanding of the global dynamics of the flow generated by the vector field. The invariant manifolds divide the phase space into regions with different behaviour. The simplest picture is in the plane, where the invariant manifolds, the separatrices, of the saddles of the system, together with the limit cycles, decompose the phase plane into connected components where the trajectories have similar $\alpha$ and $\omega$ limit sets. A standard reference on invariant manifold theory is \cite{HPS77}. In higher dimensions the structure of the invariant manifolds is, typically, much more complicated. The fundamental theorem about hyperbolic saddles is the stable (unstable) manifold theorem, see e.g. \cite{HPS77,PdW82}, which states that locally at a hyperbolic fixed point there exist manifolds of dimensions $d_s$ and $d_u$, denoting the number of negative and positive eigenvalues of the linearisation of the vector field at the fixed point, such that the tangent spaces of the stable and unstable manifolds at the fixed point are the negative and positive eigenspaces, respectively, of the linearisation. In addition, these manifolds can locally be described as graphs of functions from the negative eigenspace to the positive eigenspace, and vice versa.

To compute global invariant manifolds, one typically starts with an approximation of the local invariant manifolds lying in the corresponding eigenspace of the linearisation, and expand the global invariant manifolds step by step from the local one. For a review of a plethora of such methods see \cite{Kea05}.

For obtaining \textit{approximations} of (local) invariant manifolds there exists many references, e.g. \cite{B90,CFL03,CFL05,MH99,S90}. Few methods exist, however, that can rigorously compute \textit{enclosures} of the local invariant manifolds, which is our current objective. Some such methods are \cite{DH97,NR93,O95,Z09}.

We present a method to compute the convergence radii together with explicit error estimates of the parametrisations of the invariant manifolds; for some methods to compute such parametrisations see e.g. \cite{CFL03,CFL05,S90}. The parametrisations that we study are constructed such that the negative and positive eigenspaces of the linearisation at the fixed point are invariants of the flow of the vector field. This is a much weaker requirement than to completely linearise the vector field, as can be done, according to Siegel's theorem \cite{S52}, under certain Diophantine conditions on the eigenvalues. The idea to compute a close to identity transformation that removes all terms necessary for the transformed equation to have this property has appeared in \cite{T02}; in \cite{T04} the resulting vector field, after this close to identity transformation, was named a robust normal form. These linearisations can be seen as a special case of the parametrisations of invariant manifolds in \cite{CFL03,CFL05}, where higher order conjugacies are also considered. The constructive method to compute convergence radii and error estimates presented in this paper, however, are much easier to implement and compute than the aposteriori convergence operator from \cite{CFL03,CFL05}. Since the case of conjugacy with a linear flow on the invariant manifolds is probably the most common, the fast and simple results from this paper, that directly translate into an algorithm to compute convergence radii and error estimates, should potentially be very useful. 

This paper is organised as follows: in Section \ref{Statements} we introduce the necessary notation, recall the necessary concepts about robust normal forms from \cite{JT09,T04}, and state our main result on the existence of analytic parametrisations. In Section \ref{Proof} we prove the main theorem. The proof of the theorem is constructive and in Section \ref{algAsp} we describe an algorithm that implements the proof. Finally, in Section \ref{Examples} we calculate the local invariant manifolds in a simple planar system, similar to the discrete system studied in \cite{NR93,O95,Z09}.

%------------------------------------------------------------------------------------------------------------------------------------------------

\section{Statement of the results}\label{Statements}
Consider a vector field in $\mathbb{R}^d$ of the following form:
\begin{equation}\label{DEq}
\dot{z} = \Lambda z + F(z),
\end{equation}
with $\Lambda\in \mathcal S$, where $\mathcal S:=\{\rm diag(\lambda_{d_s}, \dots, \lambda_1, \mu_1, \dots, \mu_{d_u} ) \, : \lambda_{d_s}\leq\dots\leq\lambda_1 <0< \mu_1\leq\dots\leq\mu_{d_u}\},$ and where $F$ is an analytic function, with $F(z)=O(z^2)$.  Note that any vector field with a saddle fixed point, with distinct real eigenvalues, can (locally) be brought into this form by an affine change of variables. We decompose $z$ in the stable and unstable coordinates, $z=(x,y)\in \mathbb{R}^{d_s} \times \mathbb{R}^{d_u}$. We use $e_i=(0,\dots,0,1,0,\dots,0)$ to denote the vector in $\mathbb R^d$ with its $i$th component equal to $1$.

The structure of the parametrisation of the invariant manifolds that we are computing is based on the close to identity change of parameters associated with the robust normal forms studied in \cite{JT09,T02,T04}. In order to simplify the formulae, we use vector and multi-index notation. In this section we revise and adapt the necessary notation and results from \cite{JT09,T02}, and state our main result.

The structure of (\ref{DEq}) implies that the stable and unstable manifolds at the origin are tangent to the coordinate axes. As discussed in the introduction, we seek a parametrisation of the stable and unstable manifolds, i.e., we want to compute maps $\phi$ and $\psi$:
$$
\phi : \mathbb{R}^{d_s} \longrightarrow \mathbb{R}^d
$$
$$
\psi : \mathbb{R}^{d_u} \longrightarrow \mathbb{R}^d
$$
that are such that:
\begin{equation}\label{stablemfd}
W^s_{\loc} = \{(\xi,0)+\phi(\xi) : \xi \in U \subset \mathbb{R}^{d_s}\},
\end{equation}

\begin{equation}\label{unstablemfd}
W^u_{\loc} = \{(0,\eta)+\psi(\eta) : \eta \in V \subset \mathbb{R}^{d_u}\}.
\end{equation}

Note, this means that the stable and unstable manifolds are \textit{not} represented as graphs; we compute \textit{parametrisations} of the invariant manifolds, i.e., we allow nonlinearities in the stable coordinates of $\phi$ and the unstable coordinates of $\psi$, respectively.

We require that $\phi=O(\xi^2)$ and $\psi=O(\eta^2)$.
The maps $\phi$ and $\psi$ determine a close to identity change of coordinates in $\mathbb{R}^d$:
\begin{equation}
\Theta(\xi, \eta) = (\xi,\eta)+\phi(\xi)+\psi(\eta) .
\end{equation}

The idea of the parametrisation is that in $(\xi, \eta)$-coordinates the local stable and unstable manifolds should be given by $E_s$ ($=\mathbb{R}^{d_s}$) and $E_u$ ($=\mathbb{R}^{d_u}$), the stable and unstable tangent spaces at the fixed point. $\xi$ should be interpreted as the nominally stable, and $\eta$ as the nominally unstable coordinates. Since $\phi$ and $\psi$ do not have any constant or linear parts, the pullback of the original vector field using $\Theta$ has the following form:
\begin{equation}\label{pullbackvf}
\Theta^*(\Lambda+F)=\Lambda+G.
\end{equation}

The formal power series for the nonlinear part of the vector field in the new coordinates is: 
\begin{equation}
G=\sum_{|m|=2}^\infty g_m \zeta^m
\end{equation}
In order for the local invariant manifolds to be of the forms (\ref{stablemfd}) and (\ref{unstablemfd}) $G$ must be of order $O(\min(|\xi|,|\eta|))$.
This means that if there exists $i$, $1\leq i \leq d$, such that $g_me_i$ is a non-zero coefficient in the formal power series of $G$, then $|m_s| \geq 1$ and $|m_u|\geq 1$. 

Thus,
\begin{equation}
G\arrowvert_{E_s} \equiv 0 \quad \textrm{and} \quad G\arrowvert_{E_u} \equiv 0
\end{equation}

We call the non-negative number $|m| = |m_u| + |m_s|$ the \textit{order} of $m$, and define the set $\tilde{\mathbb N}^2 = \{m\in \mathbb N^2\colon |m|\ge 2\}$.

We split the space of multi-indices into the sets
\begin{eqnarray*}
\mathbb{V} & := & \{m\in \tilde{\mathbb N}^2 \, : |m_s| = 0 \textrm{ or } |m_u| = 0\},\\
\mathbb{U} & := & \{m\in \tilde{\mathbb N}^2 \, : |m_s|\geq 1 \textrm{ and } |m_u| \geq 1\},
\end{eqnarray*}
where $\mathbb V$ is further decomposed into $\mathbb V=\mathbb{V}_s\cup \mathbb V_u$, where 
$$
\mathbb{V}_s=\{m \in \mathbb V \, :|m_u|=0\}\quad \textrm{and} \quad \mathbb{V}_u=\{m \in \mathbb V \, :|m_s|=0\}.
$$
 
We can now define the set of admissible linear parts of (\ref{DEq}) that we consider:
$$
\mathcal F_s := \{\Lambda \in \mathcal S \, : m\in \mathbb V \Rightarrow m_s\lambda-\lambda_i\neq 0, 1\leq i \leq d_s\},
$$
$$
\mathcal F_u := \{\Lambda \in \mathcal S \, : m\in \mathbb V \Rightarrow m_u\mu-\mu_i\neq 0, 1\leq i \leq d_u\},
$$
$$
\mathcal F=\mathcal F_s \cap \mathcal F_u.
$$

We will often use the notion of \textit{filters} of a (formal) power-series: if $f(z)=\sum_{|m|\geq 2} \alpha_m z^m$, we define
$$
[f]_{\mathbb U}=\sum_{m\in \mathbb U} \alpha_m z^m,\quad [f]_{\mathbb V}=\sum_{m\in \mathbb V} \alpha_m z^m, \text{ and } [f]_{m}=\alpha_m.
$$
Note that in this notation we have
$$
[\phi]_{\mathbb V_s}=\phi, \quad [\psi]_{\mathbb V_u}=\psi, \quad \textrm{and} \quad [G]_{\mathbb U}=G.
$$

Also, we let $f^{[d]}$ denote the partial sum of the terms of $f$ up to order $d$. We use the norms $|z|=\max_{1\leq i\leq d}{\{|z_i|\}}$ and $||f||_r = \max\{|f(z)| \, : |z|<r\}$. The $r$-disc is denoted by $\mathfrak B_r$. If $X$ is a set and $r\in \mathbb{R}$, we denote by $rX$ the set $\{rx : x\in X\}$. The smallest integer, $n$, larger than a real number, $r$, is denoted by $n=\lceil r\rceil$. 

Let $\alpha=(\alpha_1, \dots, \alpha_k)\in\mathbb{R}^k$, we say that $\alpha$ is \textit{A-finitely rationally independent} if $\alpha_i\neq m\alpha$, for all multi-indices $m$ such that $2\leq |m| \leq A$, and all $1\leq i\leq k$. Finally, let $\Omega : \mathbb{Z}^{+} \rightarrow \mathbb{R}_{\geq 0}$ be defined as
\begin{equation}\label{Omegadef}
\Omega(k):=\min\left(\left|k\lambda_1-\lambda_{d_s}\right|,\left|k\mu_1-\mu_{d_u}\right|\right).
\end{equation}
We will use the following lemma, which essentially is a reformulation of \cite[Lemma 5.1]{T02}.
\begin{lemma}\label{gEst}
Assume that $\lambda$ is $\left\lceil\frac{\lambda_{d_s}}{\lambda_1}\right\rceil$-finitely rationally independent and $\mu$ is $\left\lceil\frac{\mu_{d_u}}{\mu_1}\right\rceil$-finitely rationally independent. Then, $\Lambda\in\mathcal F.$ Furthermore, for all multi-indices $m\in\mathbb{V}$ with orders
$|m|\geq \max\left(\left\lceil\frac{\lambda_{d_s}}{\lambda_1}\right\rceil,\left\lceil\frac{\mu_{d_u}}{\mu_1}\right\rceil\right)$ we have the following sharp lower bound:
\begin{equation}
|m\cdot(\lambda,\mu)-\nu|\geq \Omega(|m|),\, \quad \textrm{{\rm for }{\rm all }} \quad \nu\in\{\lambda_i\}\cup\{\mu_i\}.
\end{equation}
\end{lemma}
From this lemma it clearly follows that $\mathcal F$ is open. In addition $\mathcal F$ has full Lebesgue measure in $\mathcal S$, since it is constructed by removing countably many lines from $\mathcal{S}$.

We are now ready to state our main theorem :
\begin{theorem}\label{thm1}
Given a system $\dot{z} = \Lambda z + F(z)$, where $F(z) = \sum_{|m|\geq2}c_mz^m$ is an analytic function, $\Lambda\in\mathcal F$, and a natural number $n_1\geq\max\left(\left\lceil\frac{\lambda_{d_s}}{\lambda_1}\right\rceil,\left\lceil\frac{\mu_{d_u}}{\mu_1}\right\rceil\right)$, there exists analytic parametrisations of the stable and unstable manifolds of the forms (\ref{stablemfd}) and (\ref{unstablemfd}), converging on the disk $\mathfrak B_{r_\Theta}$, with 
\begin{equation}
\phi(\xi)\in\sum_{2\leq |m_s| \leq n_1} \alpha_{m_s} \xi^{m_s} + r_\Theta\left(\frac{|\xi|}{r_\Theta}\right)^{n_1+1}\left(1-\frac{|\xi|}{r_\Theta}\right)^{-1}\times\mathfrak{B}_1
\end{equation}
\begin{equation}
\psi(\eta)\in\sum_{2\leq |m_u| \leq n_1} \beta_{m_u} \eta^{m_u} + r_\Theta\left(\frac{|\eta|}{r_\Theta}\right)^{n_1+1}\left(1-\frac{|\eta|}{r_\Theta}\right)^{-1}\times\mathfrak{B}_1
\end{equation}
for a computable positive real number $r_\Theta$.
\end{theorem}

To prove the convergence of the change of variables $\Theta$ we proceed as in e.g. \cite{H76,SM71}, and use the method of majorants. If $f,g :\mathbb {C}^d\rightarrow \mathbb {C}^d$ are two formal power series, and $|f_m|<g_m$ for all multi-indices $m$, and all the coefficients of $g$ are real and positive, we say that $g$ {\it majorises} $f$, denoted by $f\prec g$. Thus, the convergence radius of $f$ is at least as large as that of $g$. We will majorise in two steps; given some $f:\mathbb {C}^d\rightarrow \mathbb {C}^d$, we construct $g:\mathbb {C}^d\rightarrow \mathbb {C}$ such that $f_i\prec g$, for all $i$, and then construct $h:\mathbb {C}\rightarrow \mathbb {C}$ such that $g(z,...,z)\prec h(z)$.

%------------------------------------------------------------------------------------------------------------------------------------------------

\section{Proof of the main theorem}\label{Proof}
Let $z=(x,y)$ be the original coordinates, and $\zeta=(\xi,\eta)$ the coordinates in the domain of $\Theta$. Recall, $\Theta(\xi,\eta)=(\xi,\eta)+\phi(\xi)+\psi(\eta),$ where $[\phi]_{\mathbb{V}_s}=\phi,$ and  $[\psi]_{\mathbb{V}_u}=\psi$. By inserting $z=\Theta(\zeta)$ into (\ref{DEq}), differentiating, and comparing the sides, we get:
$$D\Theta(\zeta) \dot \zeta=\dot z =\Lambda\Theta(\zeta)+F(\Theta(\zeta)).$$
Inserting this expression into (\ref{pullbackvf}) yields: 
$$D\Theta(\zeta)\Lambda \zeta+D\Theta(\zeta)G(\zeta)=\Lambda\Theta(\zeta)+F(\Theta(\zeta)),$$
we reorder the terms and get:
\begin{equation}\label{rec}
D(\phi(\xi)+\psi(\eta))\Lambda \zeta-\Lambda(\phi (\xi)+\psi(\eta))=F(\Theta(\zeta))-D\Theta(\zeta)G(\zeta).
\end{equation}
Let $L_{\Lambda}$ and $K_{\Lambda}$ be the operators 
\begin{equation}\label{Lop}
L_{\Lambda}\phi=[D\phi(\xi)\Lambda \zeta-\Lambda\phi(\xi)]_{\mathbb{V}_s},
\end{equation}
\begin{equation}\label{Kop}
K_{\Lambda}\psi=[D\psi(\eta)\Lambda \zeta-\Lambda\psi(\eta)]_{\mathbb{V}_u},
\end{equation}
where we note that 
\begin{equation}\label{Ldiv}
L_{\Lambda}(\xi^{m_s}e_i)=(m_s\lambda-(\lambda,\mu)_i)\xi^{m_s}e_i
\end{equation} 
\begin{equation}\label{Kdiv}
K_{\Lambda}(\eta^{m_u}e_i)=(m_u\mu-(\lambda,\mu)_i)\eta^{m_u}e_i.
\end{equation}
Thus, since $\phi$ is a series in $\xi^{m_s}$ terms, and $\psi$ is a series in $\eta^{m_u}$ terms, the left hand side of (\ref{rec}) can be written as 
$L_{\Lambda}\phi+K_{\Lambda}\psi$. Furthermore,
$$[L_{\Lambda}\phi+K_{\Lambda}\psi]_\mathbb{U}\equiv 0.$$ 
This means that we have to construct $G$ such that 
\begin{equation}
[F(\Theta(\zeta))-D\Theta(\zeta)G(\zeta)]_\mathbb{U}\equiv 0.
\end{equation}
Recall, we want to compute a normal form (\ref{pullbackvf}) such that $G=O(\min(|\xi|,|\eta|))$, i.e., $[G]_\mathbb{U}=G$. To be able to do this we have to construct $\phi$ and $\psi$ such that 
\begin{equation}\label{LK}
L_{\Lambda}\phi+K_{\Lambda}\psi=[F(\Theta(\zeta))-D\Theta(\zeta)G(\zeta)]_\mathbb{V}=[F(\Theta(\zeta))]_\mathbb{V}.
\end{equation}
%The terms in (\ref{DEq}), which we want to cancel with $\Theta$, come from $F$. 
To be able to simplify (\ref{LK}) by decoupling the various terms into groups, we note that  $\mathbb V_s$ and $\mathbb V_u$ are invariant under $F$, i.e., 
if $m\in\mathbb V_s$, then $m=(m_s,0)$, and for any $i$ 
$$F(z^me_i)=F(x^{m_s}e_i)=\sum_{|n|=2}^\infty c_n(x^{m_s}e_i)^n,$$
thus
$$[F(x^{m_s}e_i)]_{\mathbb V_s}=F(x^{m_s}e_i),
$$
and similarly for $\mathbb V_u$. Therefore, by filtering on the component level, we get the following two functional equations for $\phi_i$ and $\psi_i$:
\begin{equation}\label{phiRec}
(L_{\Lambda}\phi)_i=[F_i((\xi,0)+\phi(\xi))]_{\mathbb{V}_s} 
\end{equation}

\begin{equation}\label{psiRec}
(K_{\Lambda}\psi)_i=[F_i((0,\eta)+\psi(\eta))]_{\mathbb{V}_u}.
\end{equation}

It follows that $G$ should solve 
\begin{equation}\label{Grec}
G(\zeta)=[F_i(\Theta(\zeta))]_{\mathbb{U}}-D(\phi(\xi)+\psi(\eta))G(\zeta). 
\end{equation}
Since $\phi$ and $\psi$ do not contain any linear terms, this means that the coefficients of the formal power series of $G$ can be computed recursively as 
$$
g_m=\left[[F_i(\Theta(\zeta))]_{\mathbb{U}}-D(\phi(\xi)+\psi(\eta))G^{[|m|-1]}(\zeta)\right]_m.
$$

Note that (\ref{phiRec}) and (\ref{psiRec}) are the same formulae as the linear case of formulae \cite[Equations (3.5)--(3.8)]{CFL03}, but we have included their derivation for completeness, and to set the notation. 

Since $\Lambda \in \mathcal F$, and $[\phi]_{\mathbb{V}_s}=\phi$ and $[\psi]_{\mathbb{V}_u}=\psi$ by construction, we can solve (\ref{phiRec}) and (\ref{psiRec}) recursively. By computing $\phi$, $\psi$, and $G$ using the recursive formulae (\ref{phiRec}),  (\ref{psiRec}), and (\ref{Grec}), we get formal power series with the properties described in the introduction, i.e.,
\begin{equation}
\Theta^*(\Lambda+F)\arrowvert_{E_s}=\textrm{diag}(\lambda_{d_s}, \dots, \lambda_1, 0, \dots, 0)(\xi,0),
\end{equation}
and
\begin{equation}
\Theta^*(\Lambda+F)\arrowvert_{E_u}=\textrm{diag}(0, \dots, 0, \mu_1, \dots, \mu_{d_u} )(0,\eta).
\end{equation}

Thus, the stable and unstable manifolds are given by the parmetrisations (\ref{stablemfd}) and (\ref{unstablemfd}), respectively.

To bound the solutions of (\ref{phiRec}) and (\ref{psiRec}) we proceed as in \cite{JT09,T02,T04}, and prove the convergence of the change of variables using majorants and induction. 
Let 
$$N=\max\left(\left\lceil\frac{\lambda_{d_s}}{\lambda_1}\right\rceil,\left\lceil\frac{\mu_{d_u}}{\mu_1}\right\rceil\right),$$ 
be the constant from Lemma \ref{gEst} from which the explicit lower bound holds. Recall that we have assumed that $n_1\geq N$.

Let 
$$
\phi_i(\xi)=\sum_{|m_s|=2}^\infty \alpha_{i,m_s}\xi^{m_s} \quad \textrm{and} \quad \psi_i(\eta)=\sum_{|m_u|=2}^\infty \beta_{i,m_u}\eta^{m_u}
$$
be the sought change of variables, where the  $\alpha_{i,m_s}$ and $\beta_{i,m_u}$ with $|m_s|,|m_u|\leq n_1$ can be computed with any method that solves (\ref{phiRec}) and (\ref{psiRec}). To majorise the functions $\phi$ and $\psi$ we construct two one-dimensional functions $\hat \phi$ and $\hat \psi$. 
Put 
$$\hat \alpha_{k}=\sum_{|m|=k} \max_{1\leq i\leq d} \{|\alpha_{i,m}|\} \quad \textrm{and} \quad \hat \beta_{k}=\sum_{|m|=k} \max_{1\leq i\leq d} \{|\beta_{i,m}|\}.$$ 
We then define 
$$\hat \phi(\omega)= \sum_{k=2}^\infty \hat \alpha_k \omega^k \quad \textrm{and} \quad 
\hat \psi(\omega)= \sum_{k=2}^\infty \hat \beta_k \omega^k.
$$

The $\hat \phi$ and $\hat \psi$ are majorants of $\phi$ and $\psi$, respectively. Although the convergence of the parametrisations of the local stable and unstable manifolds can be proved separately, for simplicity of the exposition, we henceforth study their convergence simultaneously. 
Therefore, let 
$$\gamma_k=\hat \alpha_k + \hat \beta_k,$$ and define the joint majorant 
$$
\chi=\sum_{k=2}^\infty \gamma_k \omega^k=\hat \phi(\omega)+\hat \psi(\omega).
$$

To calculate $\alpha_{i,m}$ and $\beta_{i,m}$, with $|m|=k$, we use the operators $L_\Lambda$ and $K_\Lambda$ defined by (\ref{Lop}) and (\ref{Kop}), respectively. Their evaluation reduces by (\ref{phiRec}) and (\ref{psiRec}) to the evaluation of $k$-Taylor models of $F_i\left((\xi,0)+\phi(\xi)\right)$ and $F_i\left((0, \eta)+\psi(\eta)\right)$, respectively. The action of $L_\Lambda$ and $K_\Lambda$ on monomials are given by (\ref{Ldiv}) and (\ref{Kdiv}), and yield the following formulae for $\alpha_{i,m}$ and $\beta_{i,m}$, respectively:
\begin{equation}\label{compAlpha}
\alpha_{i,m_s} = \frac{\left[F_i\left((\xi,0)+\phi^{[k-1]}(\xi)\right)\right]_{m_s}}{\lambda m_s - (\lambda,\mu)_i},
\end{equation}

\begin{equation}\label{compBeta}
\beta_{i,m_u} = \frac{\left[F_i\left((0,\eta)+\psi^{[k-1]}(x)\right)\right]_{m_u}}{\mu m_u - (\lambda,\mu)_i}.
\end{equation}

Note that the coefficients at a certain level only depend on the previous levels. The reason is that $F$ does not contain constant or linear terms. This in turn allows for a recursive solution scheme of (\ref{phiRec}) and (\ref{psiRec}), given by (\ref{compAlpha}) and (\ref{compBeta}), respectively.

If $n_1$ is sufficiently large, then the first $n_1$ terms of $\phi$ and $\psi$ produce a good approximation of a majorant $\chi$, and we use this to determine an approximate radius of convergence for $\chi$. The validity of this radius of convergence will now be proved.
As a first step we determine, using a least squares estimator, constants $C$ and $M$ such that 
\begin{equation}\label{leastSq}
\gamma_k\leq CM^k,\quad \left\lfloor \frac{n_1}{2} \right\rfloor<k\leq n_1.
\end{equation}
The reason why we only use the terms from $\left\lfloor \frac{n_1}{2} \right\rfloor$ and onwards, is that $\left\lfloor \frac{n_1}{2} \right\rfloor$ should be large enough to capture transient phenomena in the sizes of the coefficients of $\chi$, so that the estimate from (\ref{leastSq}) is a tight bound on the coefficients in the tail of the power series of $\chi$.
The least squares estimation is done in two steps: first a standard least squares approximation is computed, then we assume that $M$ has been well approximated and increase C until (\ref{leastSq}) holds. Thus, a candidate radius of convergence is 
\begin{equation}
r_\Theta=\frac{1}{M},
\end{equation}
which needs to be verified. 

We will consider a slightly larger majorant of $\chi$. 
If $$F(z)=\sum_{|m|=2}^\infty c_mz^m\ ,$$ 
we define $$\hat c_{k}:=\sum_{|m|=k} \max_{1\leq i \leq d} \left\{|c_{i,m}|\right\},$$ and set 
$$\hat F(\omega)=\sum_{k=2}^\infty \hat c_k \omega^k.$$ $\hat F$ is clearly a majorant of $F_i(z,\dots,z)$. 
Recall the definition of (\ref{Omegadef}), and let
\begin{equation}
\Omega=\min{\left(\min_{2\leq |m| < N, \nu\in\{\lambda_i\}\cup\{\mu_i\}.}{\frac{|m\cdot(\lambda,\mu)-\nu|}{|m|}}, \frac{\Omega(N)}{N}\right)}.
\end{equation}
Note that $\frac{\Omega(k)}{k}$ is monotonically increasing for $k\geq N$, and $|m\cdot(\lambda,\mu)-\nu| \geq \Omega(k) $ 
for all $|m|\geq N$ and $\nu\in\{\lambda_i\}\cup\{\mu_i\}$.
Hence, $$|m\cdot(\lambda,\mu)-\nu|\geq \Omega |m|,$$ for all $|m|\geq 2$, and all $\nu\in\{\lambda_i\}\cup\{\mu_i\}$.
We recursively define a majorant $$\sigma(\omega)=\sum_{k=2}^\infty \delta_k \omega^k$$ of $\chi(\omega)$ by:
$$
\delta_k=\frac{1}{\Omega k}\left[\hat F(\omega+\sigma^{[k-1]}(\omega))\right]_k, \quad \textrm{for} \quad k\geq 2.
$$

In our proof of convergence of the parametrisations we will use a quadratic bound on $\hat F$. If the convergence radius of $F$ is $s$, we choose two other radii $0<s''<s'<s$, and use Cauchy-type estimates on the $\rho$-tail of $F$ on $\mathfrak{B}_{s'}$ valid on $\mathfrak{B}_{s''}$, and then require that $2r_\Theta\leq s''$. 

Indeed, let $N_d(k)$ denote the number of $d$-dimensional multiindices with absolute value $k$, then for any $i$ and $m$, 
$$|c_{i,m}| \leq \frac{ \|f\|_{s'}}{(s')^{|m|}},$$ 
and hence 
$$\hat c_k \leq N_d(k)\frac{ \|f\|_{s'}}{(s')^{k}}.$$
Therefore, we can bound the $\rho$-tail of $F$ on $s''$ as follows, 
\begin{equation}
\begin{array}{ccl}
\left|\sum_{k=\rho+1}^\infty \hat c_k \omega^k\right| & \leq & 
||F||_{s'}\sum_{k=\rho+1}^\infty N_d(k) \left(\frac{\omega}{s'}\right)^k \\
& \leq & \left(\frac{||F||_{s'}}{(s')^2}\sum_{k=\rho+1}^\infty N_d(k) \left(\frac{s''}{s'}\right)^{k-2} \right)\omega^2.
\end{array}
\end{equation} 

We denote the bound on the tail 
$$A_{s''}^\rho=\frac{||F||_{s'}}{(s')^2}\sum_{k=\rho+1}^\infty N_d(k) \left(\frac{s''}{s'}\right)^{k-2},$$ 
and define:

\begin{equation}\label{Adef}
A_{s''}=\sum_{k=2}^\rho \hat c_{k} (s'')^{k-2} + A_{s''}^\rho.
\end{equation}
Clearly, $|\hat F(z)|\leq A_{s''}|z|^2$, on $\mathfrak{B}_{s''}$.

If possible put $s''=2 r_\Theta$, otherwise take $s''$ as large as possible and put $r_\Theta=s''/2$.
If $r_\Theta > \frac{\Omega}{4 A_{2 r_\Theta}}$, decrease $r_\Theta$ until 
\begin{equation}\label{rTreq}
r_\Theta \leq \frac{\Omega}{4A_{2 r_\Theta}}.
\end{equation}

To prove the convergence of $\sigma$ we proceed as in \cite{H76,T04}. The definition of the $\delta_k$'s imply that (formally) the following equation holds:

$$\sum_{k=2}^\infty \Omega k \delta_k \omega^k = \sum_{k=2}^\infty \hat c_k\left(\omega + \sigma(\omega)\right)^k,$$
where we note that the left hand side is equal to $\Omega \omega \sigma'(\omega)$. Hence, $\sigma$ satisfies the following differential equation

$$\sigma'(\omega)=\frac{\hat F(\omega + \sigma(\omega))}{\Omega \omega},\quad \sigma(0)=0.$$

The fact that neither $\sigma$ nor $\hat F$ have any constant or linear parts imply that the following inequalities holds for any partial sum (the first is an inequality since the right hand side, in general, also includes some higher order terms):

$$0\leq \sigma'^{[k]}(\omega) \leq \frac{\hat F(\omega + \sigma^{[k-1]}(\omega))}{\Omega \omega},\quad 0\leq \omega,$$
and
$$0\leq\sigma^{[k]}(\omega) \leq \omega \sigma'^{[k]}(\omega),\quad 0\leq \omega.$$
Together they imply that
$$0\leq\sigma^{[k]}(\omega)\leq\frac{\hat F(\omega + \sigma^{[k-1]}(\omega))}{\Omega},\quad 0\leq \omega.$$

We will now use our quadratic bound on $\hat F$. Assume that $\sigma^{[k-1]}(r_\Theta)\leq r_\Theta$, for some $k$ (this trivially holds for $k=2$), then

\begin{eqnarray*}
0&\leq& \sigma^{[k]}(r_\Theta) \\ 
&\leq& \frac{\hat F(r_\Theta + \sigma^{[k-1]}(r_\Theta))}{\Omega} \\
&\leq &  \frac{\hat F(r_\Theta + r_\Theta)}{\Omega} \\
&\leq & \frac{A_{2r_\Theta}(r_\Theta + r_\Theta)^2}{\Omega} \leq r_\Theta,
\end{eqnarray*}
where the last inequality is due to (\ref{rTreq}). Hence, by induction, $\sigma(r_\Theta)\leq r_\Theta$, and since all the coefficients $\delta_k$ are positive this implies that $\sigma$ is analytic on $\mathfrak{B}_{r_\Theta}.$ It follows from the convergence of $\sigma$, by tracing the sequence of majorisations backwards, that $\Theta$ is analytic on $\mathfrak{B}_{r_\Theta}$.

The remainder terms in Theorem \ref{thm1} are found by using Cauchy bounds on $\sigma$. Since $\sigma(\omega)\leq r_\Theta$ on $\mathfrak{B}_{r_\Theta}$, the convergence radius of $\sigma$ is larger than $r_\Theta$, and we have that
$$\delta_k = \frac{1}{2\pi i}\int_{|\omega|=r_\Theta} \frac{\sigma(\omega)}{\omega^{k+1}}\, d\omega \leq r_\Theta r_\Theta^{-k}.$$
Since we use the supremum norm, the uncertainties can appear in any component simultaneously. Therefore, the remainders are added as $\mathfrak{B}_1$ scaled with geometric bounds given by the bound on the growth of the $\delta_k$'s.

%------------------------------------------------------------------------------------------------------------------------------------------------

\section{Algorithmic aspects}\label{algAsp}
The main application of the convergence proof given in this paper is that it is constructive and suitable for implementation on a digital computer. We summarise the key points of the proof given in the last section and compile it into an algorithm, computing $r_\Theta$ from the formulation of Theorem \ref{thm1}. Since $r_\Theta$ tends to be slightly pessimistic, the algorithm also computes the constants $C$ and $M$. They are candidates for the geometric bound on the tail of $\Theta$. In general $r_\Theta < \frac{1}{M}$. This algorithm is given as Algorithm \ref{mainAlgorithm}.
\begin{small}
\begin{algorithm}[ph]\label{mainAlgorithm}
 \KwData{$\Lambda$, $F$, $n_1$}
\KwResult{$\phi^{[n_1]}$, $\psi^{[n_1]}$, $r_{\Theta}$, $C$, $M$}
\For{$k=2$ \KwTo $n_1$}{
\For{$i=1$ \KwTo $d$}{
\ForAll{$|m_s|=k$}{
$\alpha_{i,m_s} = \frac{\left[F_i\left((\xi,0)+\phi^{[k-1]}(\xi)\right)\right]_{m_s}}{\lambda m_s - (\lambda,\mu)_i}
$
}
\ForAll{$|m_u|=k$}{
$\beta_{i,m_u} = \frac{\left[F_i\left((0,\eta)+\psi^{[k-1]}(x)\right)\right]_{m_u}}{\mu m_u - (\lambda,\mu)_i}
$
}
}
$\gamma_{k}=\sum_{|m_s|=k} \max_{1\leq i\leq d} \{|\alpha_{i,m_s}|\} +\sum_{|m_u|=k} \max_{1\leq i\leq d} \{|\beta_{i,m_u}|\}.$ 
}
$n_0=\lfloor n_1/2\rfloor$\;
\For{$k=n_0+1$ \KwTo $n_1$}{
$B(k,1)=1$, $B(k,2)=k$, $b(k)=\log \gamma_k$
}

$(\log C, \log M)=(B^TB)^{-1}b$ \;
\For{$k=n_0+1$ \KwTo $n_1$}{
\If{$\gamma_k>CM^k$}{
$C=\gamma_k/M^ k$
}
}
\For{$k=2$ \KwTo $\rho$}{
$\hat c_{k}=\sum_{|m|=k} \max_{1\leq i \leq d} \left\{|c_{i,m}|\right\}$
}
$\Omega=\min{\left(\min_{2\leq |m| < N, \nu\in\{\lambda_i\}\cup\{\mu_i\}.}{\frac{|m\cdot(\lambda,\mu)-\nu|}{|m|}}, \frac{\Omega(N)}{N}\right)}$\;
$r_\Theta=\frac{1}{M}$ \;
\Repeat{$Converges=True$}{
Compute $A^\rho_{2r_\Theta}$\;
$A_{2r_\Theta}=\sum_{k=2}^\rho \hat c_{k} (2r_\Theta)^{k-2} + \frac{A_{2r_\Theta}^\rho}{(2r_\Theta)^2}$\;
\eIf{$r_\Theta \leq \frac{\Omega}{4A_{2r_\Theta}}$}
{$Converges=True$}
{$Converges=False$\;
$r_\Theta=0.95r_\Theta$\;
}

}
\caption{Implementation of the proof of the main Theorem} 
\end{algorithm}
\end{small}

%------------------------------------------------------------------------------------------------------------------------------------------------

\section{Example}\label{Examples}\label{ex1}
An algorithm proving the conditions of Lemma \ref{gEst} and Theorem \ref{thm1} has been implemented in a {\tt C++} program using the C-XSC package \cite{CXSC,HH95} for interval arithmetic \cite{AH83,Mo66,Mo79,Ne90}. For automatic differentiation \cite{G00} we use a modified version of the Taylor arithmetic package \cite{BHK05}.

There are several methods to compute local invariant manifolds of discrete dynamical systems, see e.g. \cite{NR93,O95,Z09}. In principle, these methods can also be used for continuous dynamical systems by studying the time--$t$ map of the flow for some $t$. To compare with the method in \cite{Z09}, which is also able to treat flows, we study a vector field of the same form as the discrete dynamical system studied in \cite{NR93,O95,Z09}.

\begin{equation}
\left(\begin{array}{c} \dot x \\ \dot y \end{array} \right) = 
\left( \begin{array}{c} 
-0.4x+x^2+y^2\\ 
1.5y-x^3+y^3 \\ 
\end{array}\right)
\end{equation}

Using $n_1=81$, we compute (the computation takes a few seconds) $r_\Theta=0.023$, $M=2.69$ and $C=0.30$. These values yield the following bound on the error terms in Theorem \ref{thm1}: 
$$0.023\left(\frac{|\zeta|}{0.023}\right)^{82}\left(1-\frac{|\zeta|}{0.023}\right)^{-1}, \quad \textrm{ for } \zeta\in \mathfrak{B}_{0.023}.$$

\begin{figure}[ht]
\begin{center}
\includegraphics[width=0.62\textwidth]{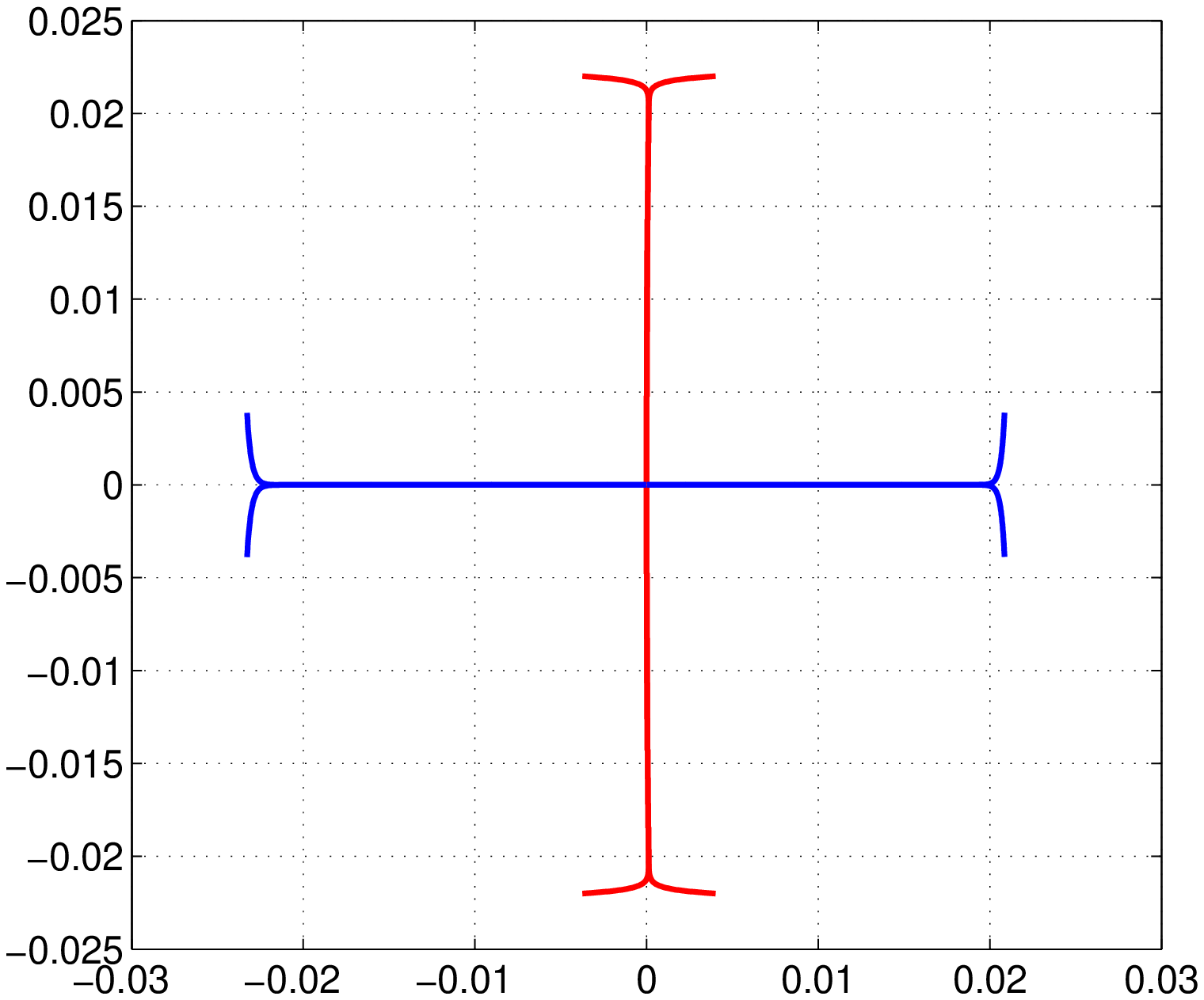}
\caption{Enclosures of the local stable (blue) and unstable (red) manifolds in the example. As the local invariant manifolds approach the converges radius of $\Theta$, the uncertainty of its location, given by the remainder term in Theorem \ref{thm1}, becomes unbounded.}\label{verified}
\end{center}
\end{figure}

\begin{figure}[h]
\begin{center}
\includegraphics[width=0.62\textwidth]{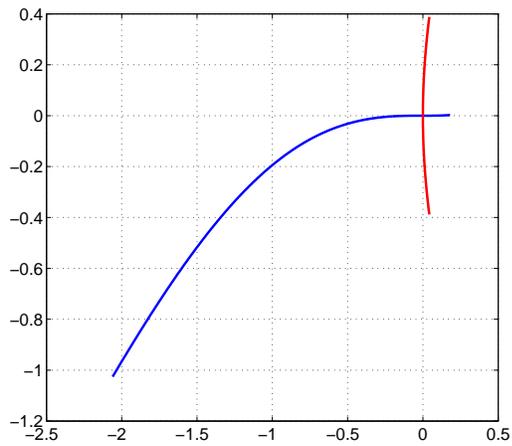}
\caption{The local stable (blue), and unstable (red) manifolds on their heuristic domain of existence.}\label{unverified}
\end{center}
\end{figure}

The image $\Theta(\mathfrak{B}_{0.023})$ is shown in Figure \ref{verified}. The image of $\Theta(\mathfrak{B}_{0.37})$ given by the candidate radius of convergence $\frac{1}{M}$ is given in Figure \ref{unverified}. By inspection we see that the image $\Theta(\mathfrak{B}_{0.023})$ contains the ball $\mathfrak{B}_{0.02}$; this can be compared with the convergence radius $0.18$ with the bound $0.241138$ on the Lipschitz constant \cite{Z08b} for the cone enclosures of the local invariant manifolds using the method from \cite{Z09}. This indicates that a method to enclose local invariant manifolds on a larger domain could be to use the method \cite{Z09} outside of the $\Theta$-image of the result of our method. 


%------------------------------------------------------------------------------------------------------------------------------------------------

\end{document}